\documentclass[12pt,reqno]{amsart}
\usepackage[T1]{fontenc}
\usepackage[utf8]{inputenc}
\usepackage{lmodern}
\usepackage{amsmath,amssymb,amsthm,mathrsfs,mathtools}
\usepackage{txfonts}
\usepackage[margin=27mm]{geometry}
\usepackage{enumitem,longtable,array,xcolor}

\ifpdf
  \usepackage[colorlinks,final,backref=page,hyperindex]{hyperref}
\else
  \usepackage[colorlinks,final,backref=page,hyperindex,hypertex]{hyperref}
\fi

\newcommand{\nc}{\newcommand}
\newcommand{\delete}[1]{}	% use this to show deleted texts {\textcolor{yellow}{#1}}

\nc{\mlabel}[1]{\label{#1}}  % Use this to suppress names
\nc{\mcite}[1]{\cite{#1}}  % Use this to suppress names
\nc{\mref}[1]{\ref{#1}}  % Use this to suppress names
\nc{\meqref}[1]{\eqref{#1}} %
\nc{\mbibitem}[1]{\bibitem{#1}} % Use this to show number name
\delete{
\nc{\mlabel}[1]{\label{#1}  % Use the next two lines to show names
{\hfill \hspace{1cm}{\small\tt{{\ }\hfill(#1)}}}}
\nc{\mcite}[1]{\cite{#1}{\small{\tt{{\ }(#1)}}}}  % Use this lines to show names
\nc{\mref}[1]{\ref{#1}{{\tt{{\ }(#1)}}}}  % Use this lines to show names
\nc{\meqref}[1]{\eqref{#1}{{\tt{{\ }(#1)}}}}  % Use this lines to show names
\nc{\mbibitem}[1]{\bibitem[\bf #1]{#1}} % Use this to show name
}

\allowdisplaybreaks[2]

\newtheorem{theorem}{Theorem}[section]
\newtheorem{proposition}[theorem]{Proposition}
\newtheorem{lemma}[theorem]{Lemma}
\newtheorem{corollary}[theorem]{Corollary}
\theoremstyle{definition}
\newtheorem{definition}[theorem]{Definition}
\newtheorem{remark}[theorem]{Remark}
\newtheorem{example}[theorem]{Example}

\newcommand{\kk}{\mathbf k}

\newcommand{\Q}{\mathbb Q}
\newcommand{\R}{\mathbb R}
\newcommand{\Z}{\mathbb Z}
\newcommand{\calf}{\mathcal F}

\newcommand{\id}{\operatorname{id}}
\newcommand{\Span}{\operatorname{span}}

\newcommand{\Ch}{\mathrm{Ch}}
\newcommand{\loc}{\top}
\newcommand{\tloc}{\otimes_{\loc}}
\newcommand{\sqrel}{\mathrel{\top_{\otimes }}}
\newcommand{\sqprod}{\odot}
\newcommand{\LX}{\mathrm L}
\newcommand{\RX}{\mathrm R}
\newcommand{\tri}{\triangleright}
\newcommand{\btri}{\blacktriangleright}
\newcommand{\eps}{\varepsilon}
\newcommand{\flip}{\mathsf T}
\newcommand{\bfi}{\mathbf i}

\font\cyrs=wncyr5

\newcommand{\shap}{{\mbox{\cyrs X}}} %sha as product
\newcommand{\rmsupp }{{\rm supp}}
\begin{document}

\title{Locality Differential Hopf Algebras and Induced Novikov Structures}

\author{Li Guo}
\address{Department of Mathematics and Computer Science, Rutgers University, Newark, NJ 07102, United States}
\email{liguo@rutgers.edu}

\author{Hongyu Xiang}
\address{School of Mathematical Sciences, Sichuan Normal University, Chengdu, 610068, China}
\email{mathxhy@aliyun.com}

\author{Bin Zhang}
\address{School of Mathematics, Sichuan University, Chengdu, 610064, China}
\email{zhangbin@scu.edu.cn}

\date{\today}

\begin{abstract}
We generalize Novikov structures to the locality setting, here locality means that algebraic operations are defined for restricted pairs of elements. Locality multi-Novikov algebras are defined and are shown to be induced from multi-derivations on locality commutative algebras. We establish a locality relation by coefficient-spaces on tensor product in order to formulate the corresponding locality co-Novikov strutures.
As a locality bialgebraic enrichment of the locality Novikov algebra, we introduce the notion of a locality bialgebra with a Novikov pair. We then show that there is a one-to-one correspondence between the set of locality biderivations on a given locality commutative bialgebra and the set of locality Novikov pairs on this locality bialgebra, which is also true for locality Hopf algebras.
The locality Hopf algebra of Chen fractions provides
examples of the construction. The partial derivations induces a biderivation, and hence a Novikov pair, on the shuffle Hopf algebra for multiple zeta values.

\end{abstract}

\subjclass[2020]{
11M32,	%Multiple Dirichlet series and zeta functions and multizeta values
16T05,	%Hopf algebras and their applications
12H05, %Differential algebra
16T30,  % Connections of Hopf algebras with combinatorics
16W25,	%Derivations, actions of Lie algebras
%17B38,	%Yang-Baxter equations and Rota-Baxter operators
16S10,	%Associative rings determined by universal properties (free algebras, coproducts, adjunction of inverses, etc.)
40B05	% Multiple sequences and series
}

\keywords{locality; differential Hopf algebras; Novikov algebra; coderivation; Chen fractions; multiple zeta values}

\maketitle

\vspace{-1cm}

\tableofcontents

\vspace{-1cm}

\section{Introduction}\mlabel{sec:introduction}
Modeled after locality differential Hopf algebras of Chen fractions, we formulate locality Hopf algebras with Novikov pairs, generalizing the classical result of Gelfand that constructs Novikov algebras from differential commutative algebras.

\subsection{Multiple zeta values and the locality algebra of Chen fractions}

Multiple zeta functions generalize the Riemann zeta function to nested sums. In the domain of convergence, they are given by
\[
 \zeta(s_1,\ldots,s_k)
 =\sum_{n_1>\cdots>n_k\ge1}
   \frac{1}{n_1^{s_1}\cdots n_k^{s_k}},
 \qquad s_1,\ldots,s_k\in\mathbb C.
\]
For positive integers $s_1,\ldots,s_k$ with $s_1\ge2$, these series converge and define the multiple zeta values (MZVs).

The study of these values goes back to Euler's work on double series
\mcite{Eu}. Since the early 1990s, MZVs and their generalizations have
been investigated in connection with number theory, algebraic geometry,
mathematical physics, quantum groups and knot theory
\mcite{BBBL,BK,Bn,GM,GZ2,H3,IKZ,Ter,Zag}.

The $\Q$-vector space spanned by MZVs is a $\Q$-subalgebra of $\R$.
The product of two MZVs admits two expansions. Their stuffle (or quasi-shuffle) product comes from multiplying their nested sums, whereas their shuffle product comes from their iterated-integral expressions, as shown by Kontsevich. Comparing the two
product formulas gives the double shuffle relations. Ihara, Kaneko and Zagier extended these to the regularized double shuffle relations \mcite{IKZ}.
The claim that the latter generate all $\Q$-linear relations among MZVs is one of the main conjectures in this area and motivates the study of other algebraic structures related to the double shuffle products.

An alternative presentation of an MZV is given by the free multiple sum:
\begin{equation}
\zeta(s_1,\ldots,s_k)
 =\sum_{n_1,\ldots,n_k\ge1}
 \frac{1}{(n_1+\cdots+n_k)^{s_1}(n_2+\cdots+n_k)^{s_2}\cdots n_k^{s_k}}.
\end{equation}
The rational functions in the summand are evaluations of
\emph{Chen fractions}. The product of two Chen fractions who sets of variables are disjoint observe the shuffle product, giving rise to the so-called locality algebra. This shuffle product leads to another proof of the shuffle relation of two MZVs \mcite{GX}.
This observation is the starting point for the locality constructions below.

\subsection{Locality differential Hopf algebra of Chen fractions}
The space of Chen fractions carries additional algebraic structures.
Partial differentiation gives a multidifferential algebra. On the other hand,
the Chen coproduct enriches the locality algebra to a locality Hopf algebra for which the partial derivatives are also coderivations \mcite{GHXZ}. Forgetting the
variable labels gives the shuffle Hopf algebra of MZVs \cite{GPZGalois}. We study the Novikov structures induced by these derivations.

Novikov algebras arose in the work of Gelfand and Dorfman on
Hamiltonian operators \mcite{GD} and of Balinsky and Novikov on
Poisson brackets of hydrodynamic type \mcite{BN85}. Let $A$ be a
commutative associative algebra with a derivation $D$. Then the
product
$$x\tri_D y:=xD(y),$$
makes $A$ into a Novikov algebra, that is,
characterized by
$$
 \begin{aligned}
 (x\tri_D y)\tri_D z-x\tri_D(y\tri_D z)
  &=(y\tri_D x)\tri_D z-y\tri_D(x\tri_D z),\\
 (x\tri_D y)\tri_D z&=(x\tri_D z)\tri_D y,
 \end{aligned}\qquad x,y,z\in A.
$$

Bruned and Dotsenko extracted the multi-indices used in regularity
structures into the notion of a multi-Novikov algebra \mcite{BrD} and showed that it is induced by commutative algebras equipped with multiple commuting derivations. Related multidifferential and multi-Novikov structures are treated in \mcite{WGZ1,GWZ}.

Restricted multiplication of Chen fractions with independent variables is a typical instance for the locality setting:
algebraic operations are considered on independent pairs, rather
than on arbitrary pairs. The framework of locality algebras and
coalgebras was developed by Clavier, Guo, Paycha and Zhang
\mcite{CGPZRen}; tensor products in this setting were studied by
Clavier, Foissy, L\'opez and Paycha \mcite{CFLP}.

The natural question for us is to identify the Novikov-type structures that
can be induced from locality derivations and biderivations on a locality bialgebra and a Hopf algebra.

\subsection{Induced locality Novikov structures}

A locality relation on a vector space $A$ is simply a symmetric relation $\top$ on $A$ with certain linearity, and a pair $a, b$ with $a\top b$ are called independent. A locality associative algebra is a locality vector space in which the product is defined only for independent pairs and satisfies the associativity where defined. We similarly define the notion of a locality (multi-)Novikov algebra.
A \emph{locality derivation} is a linear map that preserves
the locality and satisfies the Leibniz rule on independent pairs.
For a commutative associative locality algebra, such a derivation
$D$ induces a locality Novikov product $a\tri_D b:=aD(b)$.
A commuting family of locality derivations similarly induces a
locality multi-Novikov structure.

The coalgebraic extension requires a tensor locality relation.
We use the coefficient spaces of tensors to define independence on
the tensor product, whereas the restricted tensor subspace
spanned by independent factors is the range of a locality
coproduct. These two constructions serve different purposes and
are established in Sections~\mref{sec:tensor-locality} and
\mref{sec:locality-coNovikov}.

Our bialgebra enrichment of the locality Novikov algebra obtained here is the \emph{Novikov pair}, consisting of two locality products, one associative and one Novikov, and two locality coproduct satisfying certain compatible conditions. A commutative locality bialgebra $B$ with a locality biderivation $D$ induces a  Novikov pair for which the original product and coproduct are retained, and the second product and coproduct are respectively:
\begin{equation}\mlabel{eq:intro-main-data}
 a\tri_D b:=aD(b)\quad(a\top b),\qquad
 \rho_D:=(\id\otimes D)\Delta.
\end{equation}

We further show that a locality bialgebra with a Novikov pair also gives rise to a locality bialgebra with a locality biderivation. In fact, there is a one-to-one correspondence between these two structures on a locality bialgebra.  This correspondence extends to locality Hopf algebras.

\subsection{Outline of the paper}
Here is the structure of the paper.

Section~\mref{sec:locality-basics} recalls locality algebras and
the Gelfand construction before considering the locality multi-derivation and multi-Novikov structures. Here generalizing the Gelfand construction has a defect that can be measured by commutators of locality derivations.

To prepare the construction the locality coproduct, Section~\mref{sec:tensor-locality} develops the coefficient-space
locality relation on tensor products, establishes the
needed tensor multiplications, and obtains tensor products of locality Novikov algebra.

Section~\mref{sec:locality-coNovikov} extends the locality Gelfand construction to the locality coalgebra side, with the presence of a locality coderivation. Here we first do not assume the cocommutativity and compute the Novikov defect. We then obtain the locality Novikov coalgebra under the cocommutativity hypothesis.

We study the induced structure of a biderivation on a locality commutative bialgebra in Section~\mref{sec:locality-novikov}, and introduce commutative
locality bialgebras with Novikov pairs. We show that on top of a locality commutative bialgebra, there is a one-to-one correspondence between the set of biderivationa and the set of Novikov pairs. This correspondence extends to locality commutative Hopf algebras.

Section~\mref{sec:chen-hopf} applies the general construction to the locality Hopf algebra of Chen fractions and their partial
derivations. The the sum of the partial derivations descends to a biderivation on the shuffle Hopf algebra of MZVs. The corresponding locality Novikov pair on the locality Hopf algebra of Chen fractions descends to a Novikov pair on the shuffle Hopf algebra of MZVs.

\begin{remark}
The locality Hopf algebras with Novikov pairs studied here differ
from the Novikov bialgebras considered by Hong, Bai and Guo.
The latter are related to Lie bialgebras by affinization
\mcite{HBG} and admit constructions from differential
antisymmetric infinitesimal bialgebras \mcite{HBG2}.
Here the original locality Hopf operations are retained with
the induced Novikov pair. This part is the locality version of the paper on induced structures of differential Hopf algebras \mcite{BGXZ}.
\end{remark}

\section{Locality algebras and locality Novikov algebras }
\mlabel{sec:locality-basics}
This section provides background on locality differential algebras and locality Novikov algebras that will be the basis of the paper. After recalling the corresponding algebraic structures, we prove the locality Gelfand construction,
including the obstruction for multi-Novikov identities.

\subsection {Novikov algebras and Gelfand construction}\mlabel{subsec:ordinary-multi-Novikov}
We recall Novikov algebra and differential algebras before introducing their locality versions.

\begin{definition}\mlabel{d:nov}
A {\bf Novikov algebra} is a vector space with a bilinear product $\tri$ whose {\bf associator}
\[
 A(x,y,z)
 :=(x\tri y)\tri z
  -x\tri (y\tri  z)
\]
satisfies
\begin{align}
  A(x,y,z)&=A (y,x,z),
 \tag{N1}\mlabel{eq:N1}\\
 (x\tri y)\tri z&=(x\tri z)\tri y.
 \tag{N2}\mlabel{eq:N2}
\end{align}
The first identity is called the {\bf left symmetry}, and the second is called the {\bf right commutativity}.
\end{definition}

Note that the left symmetry property defines a (left) pre-Lie algebra. Thus if $(N,\tri)$ is a Novikov algebra, then
$[x,y]_\tri=x\tri y-y\tri x$ is a Lie bracket.

We have the Gelfand construction for Novikov algebras.
\begin{proposition} \mcite{GD}
\mlabel{t:clagel}
Let $(A,\cdot)$ be a commutative associative algebra, and let $D$ be a derivation of $A$:
$$D:A\to A, \ D(xy)=D(x)y+xD(y),$$
then $x\tri_Dy=xD(y)$ defines a Novikov algebra structure on $A$.
\end{proposition}

We also have multiple version of the Novikov algebra and the Gelfang construction.

\begin{definition}\mlabel{def:multi-novikov}
Let $\Omega$ be an index set.
A multi-Novikov algebra has bilinear products $\{\tri_\omega\}_{\omega\in\Omega}$
whose mixed associators
\[
 A_{\omega,\tau}(x,y,z):
 =(x\tri_\omega y)\tri_\tau z
  -x\tri_\omega(y\tri_\tau z)
\]
satisfy
\begin{align}
 A_{\omega,\tau}(x,y,z)&=A_{\omega,\tau}(y,x,z),
 \tag{MN1}\mlabel{eq:MN1}\\
 A_{\omega,\tau}(x,y,z)&=A_{\tau,\omega}(x,y,z),
 \tag{MN2}\mlabel{eq:MN2}\\
 (x\tri_\omega y)\tri_\tau z
 &=(x\tri_\tau z)\tri_\omega y.
 \tag{MN3}\mlabel{eq:MN3}
\end{align}
\end{definition}

\begin{proposition}
\mlabel{thm:multi-gelfand}
\mcite{BrD}
Let $A$ be a commutative associative algebra with pairwise commuting derivations $\{D_\omega\}$.
Then $x\tri_\omega y=xD_\omega(y)$ defines a multi-Novikov structure.
\end{proposition}

\subsection{Locality Novikov algebras}
\mlabel{subsec:locN}
\mlabel{sec:locality-derived}
The classical construction above extends to the locality setting. First let us recall the basic background on locality, as formulated in \mcite{CGPZRen,GPZGalois,GPZLoc}.

\begin{definition}
A {\bf locality set} consists of a set $X$ equipped with a symmetric relation $\loc\subseteq X\times X$. The locality set is also denoted by $X\times_\top X$.
If $x\loc y$, we say that $x$ and $y$ are {\bf independent}.
For subsets $U,V\subseteq X$, the notation $U\top V$ means that every element of $U$ is independent of every element of $V$.
A map $f:(X, \top _X)\to (Y, \top _Y)$ is a {\bf locality map} if $x_1\top _X x_2 $ implies $f(x_1)\top _Y f(x_2)$.

In a locality set $(X, \top )$, for $U\subseteq X$, we define its {\bf polar set} by
\[
 U^\loc=\{x\in X:x\loc u\text{ for every }u\in U\}.
\]
If $X=V$ is a vector space and every polar set $U^\loc$ is a linear subspace,
then $(V,\loc)$ is called a {\bf locality vector space}.
\end{definition}

\begin{definition}
A {\bf locality algebra} is a locality vector space $(A,\top)$ and a
partially defined product $m:A\times_\top A\to A$, denoted by $xy$,
such that, for every $U\subseteq A$,
\begin{equation}\mlabel{eq:locality-product-polar}
 x,y\in U^\top,\ x\top y\quad\Longrightarrow\quad xy\in U^\top,
\end{equation}
and, for all $z\in A$, the maps $x\mapsto xz$ and $x\mapsto zx$
are linear on $\{z\}^\top$. The locality algebra is {\bf associative} if
\begin{equation}\mlabel{eq:locality-assoc}
 (xy)z=x(yz)\qquad(x,y,z\text{ pairwise independent}),
\end{equation}
and {\bf commutative} if $xy=yx$ for $x\top y$.
A locality unit $1\in A$ satisfies
$$1\top A, \quad 1x=x1=x \quad (x\in A).$$
\end{definition}

We write
\[
 V\otimes_\top W:=\operatorname{span}\{v\otimes w\,|\,v\in V,\ w\in W,\ v\top w\}
 \subseteq V\otimes W
\]
for the subspace of $V\otimes W$ generated by pure tensors with independent factors, when $V,W$
are locality subspaces of a common locality vector space.

\begin{definition}\mlabel{def:locder}
Let $A$ be an associative locality algebra. A linear map $D:A\to A$
is a {\bf locality derivation} if
\begin{equation}\mlabel{eq:locality-derivation}
 \begin{gathered}
  D(U^\top)\subseteq U^\top\qquad(U\subseteq A),\\
  D(xy)=D(x)y+xD(y)\qquad(x\top y).
 \end{gathered}
\end{equation}
A family $(D_\omega)_{\omega\in\Omega}$ of locality derivations is {\bf commuting} if
$[D_\omega,D_\tau]=0$ for all $\omega,\tau\in\Omega$. The first part of this condition is also called the {\bf polar stability}.
\end{definition}

Polar stability is stronger than merely requiring $D$ to be a locality map:
it directly yields $x\loc y\Rightarrow D(x)\loc y$, and consequently, $D(x)\loc D(y)$.

\begin{remark}
Every finite composition $D_{\omega_1}\cdots D_{\omega_r}$ is polar stable.
Consequently, if $x_1,\ldots,x_n$ are pairwise independent,
then applying an arbitrary finite composition of derivations to each input preserves their pairwise independence.
\end{remark}

\begin{definition}\mlabel{def:loc-mn}
A {\bf locality multi-Novikov algebra} is a locality vector space $(N, \top)$
equipped with a family of locality bilinear operations $\{\tri_\omega\}_{\omega\in \Omega}$ on $ N\times _\top N$, satisfying for $\omega, \tau \in \Omega,$
\[
 x,y\in U^\loc,\quad x\loc y
 \quad\Longrightarrow\quad x\tri_\omega y\in U^\loc,
 \tag{LP}\mlabel{eq:lmn-polar-stability}
\]
and, for pairwise independent $x,y,z$, the mixed associators
satisfy
\begin{align}
 A_{\omega,\tau}(x,y,z)&=A_{\omega,\tau}(y,x,z),
 \tag{LMN1}\mlabel{eq:locality-NCMN1}\\
 A_{\omega,\tau}(x,y,z)&=A_{\tau,\omega}(x,y,z),
 \tag{LMN2}\mlabel{eq:locality-MN2}\\
 (x\tri_\omega y)\tri_\tau z
 &=(x\tri_\tau z)\tri_\omega y.
 \tag{LMN3}\mlabel{eq:locality-NCMN3}
\end{align}
Here condition \textup{(LP)} ensures that all these nested products are defined.

Without condition (LMN2), it is called a {\bf locality noncommuting multi-Novikov algebra}.
When there is only one index, it is called a {\bf locality Novikov algebra}.
\end{definition}

\begin{theorem}\mlabel{thm:gelfand}
Let $A$ be a locality commutative associative algebra, and let $\{D_\omega\}$ be a family of locality derivations.
Define a family of operators  $\{\tri _\omega \}$ by
$$x\tri_\omega y:=xD_\omega(y),  x\top y.$$
Then the family satisfy \textup{(LP)}, \textup{(LMN1)}, and
\textup{(LMN3)}, while the defect in \textup{(LMN2)} is
\begin{equation}
 A_{\omega,\tau}(x,y,z)-A_{\tau,\omega}(x,y,z)
 =-xy[D_\omega,D_\tau](z).
 \mlabel{eq:locality-commutator-defect}
\end{equation}
Thus pairwise commuting locality derivations yield a locality multi-Novikov structure.
Conversely, if $A$ has a locality unit and these operations satisfy \textup{(LMN2)},
then the derivations commute pairwise.
\end{theorem}

\begin{proof}
If $x\loc y$, then preservation of polars gives $x\loc D_\omega(y)$. So $\tri _\omega$ is defined on $A\times _\top A$, and locality bilinearity follows.

In addition,  if $x,y\in U^\loc$, then
$D_\omega(y)\in U^\loc\cap\{x\}^\loc$.
Polar stability of the original multiplication gives
$xD_\omega(y)\in U^\loc$, proving \textup{(LP)}.

Fix pairwise independent elements $x,y,z$, we have
\begin{align}
 (x\tri_\omega y)\tri_\tau z
 &=xD_\omega(y)D_\tau(z),\mlabel{eq:A1}\\
 x\tri_\omega(y\tri_\tau z)
 &=xD_\omega(y)D_\tau(z)+xyD_\omega D_\tau(z).
 \notag
\end{align}
Subtracting the two expressions yields
\begin{equation}
 A_{\omega,\tau}(x,y,z)=-xyD_\omega D_\tau(z).
 \mlabel{eq:assoc}
\end{equation}
So
\begin{itemize}
    \item the right-hand side is symmetric in $x,y$, proving \textup{(LMN1)};
    \item equation~\meqref{eq:A1} is invariant under simultaneously interchanging $(y,\omega)$ and $(z,\tau)$,
    which proves \textup{(LMN3)};
    \item interchanging the operator indices and subtracting gives the stated defect formula.
\end{itemize}
This proves that the commutativity of the derivations is sufficient for \textup{(LMN2)}.

For the converse, just set $x=y=1$ in the defect formula.
This substitution is legitimate because a locality unit is independent of every element.
\end{proof}

\section{Coefficient-space locality on tensor products}
\mlabel{sec:tensor-locality}
To study the coproduct in locality setting, we have to deal with locality relations in tensor products.

\subsection{The locality relation on a tensor product}
Fix locality vector spaces $(V,\top_V)$ and $(W,\top_W)$.
\begin{definition}\mlabel{def:coefficient-spaces}
For $x \in V\otimes W$, we define
\begin{equation}\mlabel{eq:coefficient-spaces}
 \LX_x:=\{(\id_V\otimes\lambda)x:\lambda\in W^*\},\qquad
 \RX_x:=\{(\varphi\otimes\id_W)x:\varphi\in V^*\}.
\end{equation}
All duals here are algebraic duals.
Then $\LX_x$ is a subspace of $V$, and $\RX_x$ is a subspace of $W$.
They are called the {\bf left} and {\bf right coefficient spaces} of $x$, respectively.
\end{definition}

Let us list properties for the coefficient spaces.
\begin{lemma}\mlabel{lem:minimal-coefficients} Let $x,y\in V\otimes W$, $P\subseteq V$, $Q\subseteq W$ be linear subspaces, and $f\in {\rm End }(V)$, $g\in {\rm End }(W)$.
\begin{enumerate}
   \item \mlabel{i:mc1} The subspaces $\LX_x$ and $\RX_x$ are finite-dimensional, and $x\in\LX_x\otimes\RX_x$.
\item \mlabel{i:mc2}  \begin{equation}\mlabel{eq:minimal-coefficients}
     x\in P\otimes Q\quad\Longleftrightarrow\quad\LX_x\subseteq P,\quad\RX_x\subseteq Q.
    \end{equation}
    \item \mlabel{i:mc3}
    \begin{gather}
    \LX_{x+y}\subseteq\LX_x+\LX_y,\quad \RX_{x+y}\subseteq\RX_x+\RX_y,\mlabel{eq:coefficient-sum}\\
    \LX_{(f\otimes g)x}\subseteq f(\LX_x),\quad
    \RX_{(f\otimes g)x}\subseteq g(\RX_x). \mlabel{eq:coefficient-map}
    \end{gather}

\end{enumerate}
\end{lemma}

\begin{proof}
\meqref{i:mc1}
Let
$$x=\sum a_i\otimes b_i.$$
Then by definition,
$$\LX _x\subseteq {\rm span}\{a_i\}, \quad\RX _x\subseteq {\rm span}\{b_i\}, $$
and so $\LX _x$, $\RX _x$ are finite dimensional.

Let $\{v_i\}$ be a basis of ${\rm span}\{a_i\} $. Then
$$x=\sum_i v_i\otimes y_i.$$
By definition,
$$\RX _x={\rm span }\{y_i\}.$$
Now taking a basis $\{w_j\}$ of $\RX _x$, then
$$x=\sum _j x_j\otimes w_j,$$
and so
$$\LX _x={\rm span}\{x_j\}.$$
Therefore,
$x$ is in $\LX _x\otimes \RX _x.$

\meqref{i:mc2} The first part of the proof of\meqref{i:mc1} shows that
$$x\in P\otimes Q \Rightarrow \LX_x\subseteq P, \RX _x\subseteq Q .$$
On the other hand, if $\LX _x\subseteq P$ and $\RX _x\subseteq Q$, then
$$x\in \LX _x\otimes \RX _x\subseteq P\otimes Q.$$
So
$$ x\in P\otimes Q\quad\Longleftrightarrow\quad\LX_x\subseteq P,\quad\RX_x\subseteq Q.$$

\meqref{i:mc3}
Now we have
$$x\in \LX _x\otimes \RX_x \subseteq (\LX_x+\LX _y)\otimes (\RX _x+\RX _y),$$
$$y\in \LX _y\otimes \RX_y \subseteq (\LX_x+\LX _y)\otimes (\RX _x+\RX _y).$$
Thus
$x+y\in  (\LX_x+\LX _y)\otimes (\RX _x+\RX _y).$
So
$ \LX_{x+y}\subseteq\LX_x+\LX_y, \RX_{x+y}\subseteq\RX_x+\RX_y.$
 For $f:V\to V$ and $g:W\to W$, since $x\in \LX_x\otimes \RX_x$, we have
$$(f\otimes g)(x)\in f(\LX_x)\otimes g(\RX_x).$$
By Eq.(\mref{eq:minimal-coefficients}), we have
$$\hspace{2.5cm}
(f\otimes g)(x)\in f(\LX_x)\otimes g(\RX_x)\Leftrightarrow \LX_{(f\otimes g)x}\subseteq f(\LX_x), \RX_{(f\otimes g)x}\subseteq g(\RX_x). \hspace{2.5cm}\qedhere$$
\end{proof}

\begin{definition}\mlabel{def:square-locality} For locality vector spaces $(V, \top _V)$ and $(W, \top _W)$, we define a locality relation on $V\otimes W$, denoted by $\top_\otimes:=(A\otimes B)\times_{\top_\otimes}(A\otimes B)$, by
\begin{equation}\mlabel{eq:tensor-relation}
 x\sqrel y\quad\Longleftrightarrow\quad
 \LX_x\top_V\LX_y,\qquad\RX_x\top_W\RX_y.
\end{equation}
\end{definition}

\begin{proposition}
\mlabel{pp:square-polar}
Relation~\meqref{eq:tensor-relation} is equivalent to the existence of finite expansions
\begin{equation}\mlabel{eq:independent-expansions}
 x=\sum_i v_i\otimes w_i,\qquad y=\sum_jv'_j\otimes w'_j,
 \qquad v_i\top_Vv'_j,\quad w_i\top_Ww'_j\quad(\forall~i,j).
\end{equation}
For $U\subseteq V\otimes W$, put
$\LX( U)=\sum_{x\in U}\LX_x$ and $\RX(U)=\sum_{x\in U}\RX_x$. Then
\begin{equation}\mlabel{eq:tensor-polar}
 U^{\top_\otimes}
 =\LX( U)^{\top_V}\otimes\RX( U)^{\top_W}.
\end{equation}
In particular, $(V\otimes W,\sqrel)$ is a locality vector space.
\end{proposition}
\begin{proof} Since $x\in \LX _x\otimes \RX_x$, $y\in \LX _y\otimes \RX_y$, and $\LX _x\top _V \LX_y$, $\RX_x\top _W \RX _y$, we can choose expansions for $x$, $y$ in $\LX_x\otimes\RX_x$ and $\LX_y\otimes\RX_y$ satisfying conditions in \meqref{eq:independent-expansions}.

Conversely, if \meqref{eq:independent-expansions} holds, then $\LX_x\subseteq\Span\{v_i\}$ and
$\LX_y\subseteq\Span\{v'_j\}$.
So
$ \Span\{v_i\} \top _V \Span\{v_j'\}$
which gives
$\LX _x\top _V\LX _y.$
Similarly we have
$\RX _x\top _W\RX _y,$
yielding the equivalence.

For an arbitrary $U$, Lemma~\mref{lem:minimal-coefficients} gives
\[
 x\in U^{\top_\otimes}
 \Longleftrightarrow
 \LX_x\subseteq\LX(U)^{\top_V},\quad
 \RX_x\subseteq\RX(U)^{\top_W}
 \Longleftrightarrow
 x\in\LX(U)^{\top_V}\otimes\RX( U)^{\top_W}.
\]
The space on the right is a linear subspace.
\end{proof}

\begin{remark}
The relation $\top_\otimes$ is imposed on the tensor product $V\otimes W$.  Its coefficient-space definition and polar formula are sufficient for the factorwise multiplication theorem below. This relation is related to the
restricted locality tensor products of Clavier--Foissy--L\'opez--Paycha
\mcite{CFLP}.
\end{remark}

\subsection{Tensor products of locality algebras}
Let $(A, \top _A, m_A)$ and $( B, \top _B, m_B)$ be locality algebras. Define the map
$$\odot:(A\otimes B)\times_{\top_\otimes}(A\otimes B)
\longrightarrow A\otimes B,$$
as follows: for $(x, y) \in \top _\otimes=(A\otimes B)\times_{\top_\otimes}(A\otimes B)$, choose independent expansions
$x=\sum_i a_i\otimes b_i$ and $y=\sum_j c_j\otimes d_j$ with $a_i\top _A c_j$, $b_i\top _B d_j$, and define
\begin{equation}\mlabel{eq:tensor-product}
 x\sqprod y:=\sum_{i,j}a_ic_j\otimes b_id_j.
\end{equation}

\begin{theorem}\mlabel{thm:tensor-algebra} The map $\odot$ is well-defined and makes $A\otimes B$ a locality algebra.
If $A$ and $B$ are unital with unit $1_A$ and $1_B$, so is $A\otimes B$ with unit $1_A\otimes1_B$.
\end{theorem}
\begin{proof}
Since $\LX_x\top_A\LX_y$, the multiplication $m_A$ is bilinear on the Cartesian product $\LX_x\times\LX_y$,
inducing
$$m_L:\LX_x\otimes\LX_y\to A,$$
Similarly we obtain
$$m_R:\RX_x\otimes\RX_y\to B.$$
Let $\sigma_{23}$ interchange the second and third tensor factors. Then we have a well-defined product on $(\LX _x \otimes \RX _x) \otimes (\LX _y\otimes \RX_y)  $,
\begin{equation}\mlabel{eq:tensor-intrinsic-product}
 (m_L\otimes m_R)\sigma_{23}:(\LX _x \otimes \RX _x) \otimes (\LX _y\otimes \RX_y) \to A\otimes B
\end{equation}

For any independent expansions $x=\sum_i a_i\otimes b_i$ and $y=\sum_j c_j\otimes d_j$ with $a_i\top _A c_j$, $b_i\top _B d_j$,
the definition gives
$$\LX _x\subseteq {\rm span}\{a_i\}, \LX _y\subseteq {\rm span}\{c_j\}, \RX _x\subseteq {\rm span}\{b_i\}, \RX _y\subseteq {\rm span}\{d_j\},$$
and similarly we have a product
$$
 (m_L\otimes m_R)\sigma_{23}:({\rm span}\{a_i\} \otimes {\rm span}\{b_i\}) \otimes ({\rm span}\{c_j\} \otimes {\rm span}\{d_j\}) \to A\otimes B.
$$
When this product is restricted to $(\LX _x \otimes \RX _x) \otimes (\LX _y\otimes \RX_y)  $, it gives \meqref{eq:tensor-intrinsic-product}, so the product is well-defined.

Now let us prove that the product $\odot$ is polar stable, that is, for $U\subseteq A\otimes B$, if $x\top y$, $x, y \in U^{\top _\otimes}$, then
$$x\odot y\in U^{\top _\otimes}.$$

Since $x \in \LX _x\otimes \RX_x$, $y\in \LX _y\otimes \RX_y$, we have expansions:
$x=\sum_i a_i\otimes b_i$ and $y=\sum_j c_j\otimes d_j$ with $a_i\in \LX _x $, $b_i\in  \RX_x$, $c_j\in \LX _y $, $d_j\in  \RX_y$
and
$$x\odot y=\sum a_ic_j\otimes b_id_j,$$
Thus we have
\begin{equation}\mlabel{eq:product-coefficients}
 \LX_{x\sqprod y}\subseteq\Span\{ac:a\in\LX_x,c\in\LX_y\} \top _A \LX(U),\quad
 \RX_{x\sqprod y}\subseteq\Span\{bd:b\in\RX_x,d\in\RX_y\}\top _B \RX (U).
\end{equation}
Thus $\odot$ is polar stable.

Locality bilinearity follows from the factorwise bilinearity of
$m_A$ and $m_B$ on the coefficient spaces. Thus the product makes
$A\otimes B$ a locality algebra. If both factor multiplications are
associative, the tensor product multiplication is also associative,
since the two triple products agree on independent pure tensors and
hence, by simultaneous coefficient-space expansions, on every
pairwise independent triple of tensors. Commutativity is inherited
when both factor multiplications are commutative.

Finally, $1_A\top_A A$ and $1_B\top_B B$, so $1_A\otimes1_B$ is independent of every tensor,
and the unit identities hold factor by factor.
\end{proof}

\subsection{Tensor products of locality algebras and locality Novikov algebras}
\mlabel{subsec:locality-tensor-novikov}

\begin{lemma}\mlabel{lem:loc-tensor-novikov}
Let $(A,\top_A,\cdot)$ be a commutative locality associative algebra
and $(B,\top_B,\tri)$ a locality Novikov algebra.
Then $A\otimes B$, equipped with the locality relation $\top_\otimes$ and the
factorwise partial product
\[
(a\otimes x)\btri(b\otimes y)
=ab\otimes(x\tri y),
\]
is a locality Novikov algebra.
\end{lemma}

\begin{proof}
Arguing as in the proof of Theorem \mref {thm:tensor-algebra}, the product is well-defined,
locality bilinear and polar-stable. So we only need to check the locality Novikov axioms.

For three tensors that are pairwise $\top_\otimes$-independent, their
left and right coefficient spaces are pairwise independent.
Choosing finite expansions in those coefficient spaces reduces both
Novikov identities to pure tensors
$X=a\otimes x$, $Y=b\otimes y$, $Z=c\otimes z$ with
pairwise independent $a,b,c$ and pairwise independent $x,y,z$.

Commutativity and associativity in $A$, together with right
commutativity in $B$, give
\[
\begin{aligned}
(X\btri Y)\btri Z
&=(ab)c\otimes((x\tri y)\tri z)\\
&=(ac)b\otimes((x\tri z)\tri y)\\
&=(X\btri Z)\btri Y.
\end{aligned}
\]
Similarly, left symmetry in $B$ gives
$$\begin{aligned}
&(X\btri Y)\btri Z-X\btri(Y\btri Z)\\
&\quad=abc\otimes
\bigl((x\tri y)\tri z-x\tri(y\tri z)\bigr)\\
&\quad=bac\otimes
\bigl((y\tri x)\tri z-y\tri(x\tri z)\bigr)\\
\hspace{4cm} &\quad=(Y\btri X)\btri Z-Y\btri(X\btri Z).
\hspace{5cm}\qedhere
\end{aligned}
$$
\end{proof}

\section{Locality coalgebras and locality co-Novikov coalgebras}
\mlabel{sec:locality-coNovikov}
\subsection{Novikov coalgebras and co-Gelfand construction}
\mlabel{subsec:novikov-co-background}
Throughout this section,
\[
 \flip(x\otimes y)=y\otimes x,\qquad
 \flip_{12}=\flip\otimes\id,
 \qquad \flip_{23}=\id\otimes\flip.
\]

\begin{definition}\mlabel{def:ordinary-Novikov-coalgebra}
A {\bf Novikov coalgebra} is a vector space $C$ equipped with a linear map
$\delta:C\to C\otimes C$
satisfying
\begin{align}
 \mathfrak a_\delta&=\flip_{12}\mathfrak a_\delta,
 \tag{NC1}\mlabel{eq:paper1-NC1}\\
 (\delta\otimes\id)\delta
 &=\flip_{23}(\delta\otimes\id)\delta.
 \tag{NC2}\mlabel{eq:paper1-NC2}
\end{align}
Here
\[
 \mathfrak a_\delta
 =(\delta\otimes\id)\delta-(\id\otimes\delta)\delta.
\]
\end{definition}

By the natural pairing $\langle, \rangle$ between $C$ and $C^*$, the operator $\delta  $ on $C$ defines an adjoint operator $\tri $ on $C^*$: for $f,g\in C^*$, put
\[
 \langle f\tri g,c\rangle=\langle f\otimes g,\delta(c)\rangle.
\]
Then \meqref{eq:paper1-NC1}--\meqref{eq:paper1-NC2} are the
transposes of \meqref{eq:N1}--\meqref{eq:N2}. In particular, for finite dimensional $C$,
\[
 (C,\delta)\text{ is a Novikov coalgebra}
 \ \Longleftrightarrow\ (C^*,\tri)\text{ is a Novikov algebra}.
\]

\begin{definition}\mlabel{def:multi-Novikov-coalgebra}
For a coassociative coalgebra $(C, \Delta)$ with a family of  linear maps $\delta_\omega:C\to C\otimes C$, $\omega\in\Omega$, set
\[
 \mathfrak a_{\omega,\tau}
 =(\delta_\omega\otimes\id)\delta_\tau
  -(\id\otimes\delta_\tau)\delta_\omega.
\]
The family defines a {\bf multi-Novikov coalgebra} if, for all
$\omega,\tau\in\Omega$,
\begin{align}
 \mathfrak a_{\omega,\tau}
 &=\flip_{12}\mathfrak a_{\omega,\tau},
 \tag{MNC1}\mlabel{eq:multi-coNov-one}\\
 \mathfrak a_{\omega,\tau}
 &=\mathfrak a_{\tau,\omega},
 \tag{MNC2}\mlabel{eq:multi-coNov-two}\\
 (\delta_\omega\otimes\id)\delta_\tau
 &=\flip_{23}(\delta_\tau\otimes\id)\delta_\omega.
 \tag{MNC3}\mlabel{eq:multi-coNov-three}
\end{align}
Omitting \textup{(MNC2)} gives a {\bf noncommuting multi-Novikov
coalgebra}. For $|\Omega|=1$, these are
\meqref{eq:paper1-NC1}--\meqref{eq:paper1-NC2}.
\end{definition}

The maps $\{\delta _\omega\}$ induce products $\{\tri _\omega\}$ on $C^*$ by the natural pairing $\langle,\rangle$. Their associators satisfy
\[
 \langle A_{\omega,\tau}(f,g,h),c\rangle
 =\langle f\otimes g\otimes h,\mathfrak a_{\omega,\tau}(c)\rangle.
\]
Thus \textup{(MNC1)}--\textup{(MNC3)} transpose the three mixed
identities in Subsection~\mref{subsec:ordinary-multi-Novikov}.

By duality, we have the co-Gelfand construction. We state it here without assuming the cocommutativity.

\begin{proposition}
\mlabel{pp:ordinary-co-Gelfand}
Let $(C,\Delta)$ be a coassociative coalgebra that is not necessarily counital or cocommutative.
Let $D_\omega$, $\omega\in\Omega$,
be coderivations:
\[
 \Delta D_\omega
 =(D_\omega\otimes\id+\id\otimes D_\omega)\Delta.
\]
Set
\[
 \delta_\omega=(\id\otimes D_\omega)\Delta,
\]
Then
\begin{align}
 \mathfrak a_{\omega,\tau}
       -\flip_{12}\mathfrak a_{\omega,\tau}
 &=-(\id\otimes\id\otimes D_\tau D_\omega)(\id-\flip_{12})\Delta^{(2)},
 \mlabel{eq:coG-cocommutativity-defect-one}\\
 \mathfrak a_{\omega,\tau}-\mathfrak a_{\tau,\omega}
 &=(\id\otimes\id\otimes[D_\omega,D_\tau])\Delta^{(2)},
 \mlabel{eq:coG-commutator-defect}\\
 (\delta_\omega\otimes\id)\delta_\tau
       -\flip_{23}(\delta_\tau\otimes\id)\delta_\omega
 &=(\id\otimes D_\omega\otimes D_\tau)(\id-\flip_{23})\Delta^{(2)}.
 \mlabel{eq:coG-cocommutativity-defect-three}
\end{align}
Here $[D_\omega,D_\tau]=D_\omega D_\tau-D_\tau D_\omega$, $ \Delta^{(2)}=(\Delta\otimes\id)\Delta=(\id\otimes\Delta)\Delta$.
\end{proposition}
\begin{proof}
The coderivation identity gives
\begin{align*}
 \delta_\tau D_\omega
 &=(\id\otimes D_\tau)
       (D_\omega\otimes\id+\id\otimes D_\omega)\Delta\\
 &=(D_\omega\otimes D_\tau+\id\otimes D_\tau D_\omega)\Delta.
\end{align*}
Thus
\begin{align*}
 (\delta_\omega\otimes\id)\delta_\tau
 &=(\id\otimes D_\omega\otimes\id)
       (\Delta\otimes\id)(\id\otimes D_\tau)\Delta\\
 &=(\id\otimes D_\omega\otimes D_\tau)\Delta^{(2)},\\
 (\id\otimes\delta_\tau)\delta_\omega
 &=(\id\otimes\delta_\tau D_\omega)\Delta\\
 &=(\id\otimes D_\omega\otimes D_\tau
       +\id\otimes\id\otimes D_\tau D_\omega)\Delta^{(2)}.
\end{align*}
Subtraction yields
$$\mathfrak a_{\omega,\tau}
 =-(\id\otimes\id\otimes D_\tau D_\omega)\Delta^{(2)}.$$

Coassociativity gives
\begin{align*}
 (\id-\flip_{12})\Delta^{(2)}
 &=((\Delta-\flip\Delta)\otimes\id)\Delta \\
 (\id-\flip_{23})\Delta^{(2)}
 &=(\id\otimes(\Delta-\flip\Delta))\Delta.
\end{align*}
Since
\begin{gather*}
 \flip_{12}(\id\otimes\id\otimes D_\tau D_\omega)
   =(\id\otimes\id\otimes D_\tau D_\omega)\flip_{12},\\
 \flip_{23}(\id\otimes D_\tau\otimes D_\omega)
   =(\id\otimes D_\omega\otimes D_\tau)\flip_{23},
\end{gather*}
These identities prove the proposition.
\end{proof}

\begin{corollary} If $C$ is cocommutative, then $(C, \Delta, \delta _\omega)$ is a noncommuting multi-Novikov
coalgebra. In this case, if the coderivations are pairwise commuting, then $(C, \Delta, \delta _\omega)$ is a  multi-Novikov coalgebra.
\end{corollary}

\begin{proof} Cocommutativity means
$$\flip \Delta =\Delta,$$
Thus we have
$$\flip _{12}\Delta ^{(2)}=\flip _{12}(\Delta \otimes \id )\Delta =(\flip \Delta \otimes \id)\Delta =\Delta ^{(2)}, $$
Similarly, we have
$$\flip _{23}\Delta ^{(2)}=\flip _{23}(\id \otimes \Delta  )\Delta =\Delta ^{(2)}, $$
Hence \textup{(MNC1)} and \textup{(MNC3)} hold.

For pairwise commuting coderivations,
$$\mathfrak a_{\omega,\tau}-\mathfrak a_{\tau,\omega}
 =(\id\otimes\id\otimes[D_\omega,D_\tau])\Delta^{(2)}=0.$$
Hence \textup{(MNC2)} holds.
\end{proof}

\begin{corollary} If $C$
has a counit, then, without assuming cocommutativity,
\[
 \mathfrak a_{\omega,\tau}=\mathfrak a_{\tau,\omega}
 \quad\Longleftrightarrow\quad [D_\omega,D_\tau]=0.
\]
\end{corollary}

\begin{proof}
For a counit $\eps$, put $E=\eps\otimes\eps\otimes\id$. Then
\[
 E\Delta^{(2)}
 =(\eps\otimes((\eps\otimes\id)\Delta))\Delta
 =(\eps\otimes\id)\Delta=\id_C,
\]
and, for every $F\in\operatorname{End}_{\kk}(C)$,
\[
 E(\id\otimes\id\otimes F)\Delta^{(2)}
 =FE\Delta^{(2)}=F.
\]
Applying $E$ to~\meqref{eq:coG-commutator-defect} gives
\[
 E(\mathfrak a_{\omega,\tau}-\mathfrak a_{\tau,\omega})
 =[D_\omega,D_\tau].
\]
This proves the counital converse.
\end{proof}

\subsection{Locality Novikov coproducts}
\mlabel{subsec:locality-coNovikov}
For a locality vector space $(C, \top)$, and subspaces $V,W\subseteq C$, set
\[
 V\tloc W=\Span\{v\otimes w:v\in V,\ w\in W,\ v\top w\},
\]
and, for $n\ge2$,
\[
 V^{\otimes_\top n}
 =\Span\{c_1\otimes\cdots\otimes c_n: c_i \in V,
                    c_i\top c_j\ (i\ne j)\}\subseteq C^{\otimes n}.
\]

\begin{definition}[{\cite[Definition~4.3]{CGPZRen}}]
\mlabel{def:locality-coalgebra}
A {\bf locality coalgebra} is a tuple $(C,\top,\Delta,\eps)$,
where $(C,\top)$ is a locality vector space and
$\Delta:C\to C\otimes C$, $\eps:C\to\kk$ are linear, satisfying
\begin{gather}
 (\Delta\otimes\id)\Delta=(\id\otimes\Delta)\Delta,
 \mlabel{eq:locality-coassociativity}\\
 (\eps\otimes\id)\Delta=\id_C=(\id\otimes\eps)\Delta,
 \mlabel{eq:locality-counit}\\
 \Delta(U^\top)\subseteq U^\top\tloc U^\top
 \qquad(U\subseteq C).
 \mlabel{eq:coalgebra-polar}
\end{gather}
Here $\kk\otimes C\cong C\cong C\otimes\kk$. The condition (\mref{eq:coalgebra-polar}) is called the {\bf polar decomposition property}.
It is {\bf cocommutative} if $\flip\Delta=\Delta$.
Omitting $\eps$ and \meqref{eq:locality-counit} gives the
noncounital version.
\end{definition}

\begin{lemma}[{\cite[Lemma~4.6(i),(iii)]{CGPZRen}}]
\mlabel{lem:coproduct-square-locality}
For $x\top y$ in a locality coalgebra, there are finite expansions
\[
 \Delta x=\sum_i a_i\otimes b_i,\qquad
 \Delta y=\sum_j c_j\otimes d_j
\]
such that
\[
 a_i\top b_i,\qquad c_j\top d_j,\qquad
 \{a_i,b_i\}_i\top\{c_j,d_j\}_j.
\]
Consequently, for the coefficient-space relation of Section~\mref{sec:tensor-locality},
\begin{equation}\mlabel{eq:coproduct-square-locality}
 \Delta(x)\sqrel\Delta(y).
\end{equation}
Moreover, for $n\ge2$ and $1\le r\le n$,
\[
 (\id^{\otimes(r-1)}\otimes\Delta\otimes\id^{\otimes(n-r)})
     (C^{\otimes_\top n})\subseteq C^{\otimes_\top(n+1)}.
\]
\end{lemma}

\begin{definition}\mlabel{def:coder}
On a locality coalgebra $(C,\top,\Delta,\eps)$, a linear map
$D:C\to C$ is a {\bf locality coderivation} if
\[
 D(U^\top)\subseteq U^\top\qquad(U\subseteq C),
\]
\begin{equation}\mlabel{eq:coder}
 \Delta D=(D\otimes\id+\id\otimes D)\Delta,
 \qquad \eps D=0.
\end{equation}
A family $(D_\omega)_{\omega\in\Omega}$ is {\bf commuting} if
$[D_\omega,D_\tau]=0$ for all $\omega,\tau$.
\end{definition}
For a counital coalgebra, the last condition in \meqref{eq:coder}
already follows from its first identity:
\[
 \eps D=(\eps\otimes\eps)\Delta D=2\eps D
 \quad\Longrightarrow\quad\eps D=0.
\]
Polar stability gives, for any finite compositions $E_i$ of the
$D_\omega$ (including the identity),
\[
 c_i\top c_j\ (i\ne j)
 \quad\Longrightarrow\quad E_i(c_i)\top E_j(c_j)\ (i\ne j).
\]

\begin{definition}
\mlabel{def:locality-Novikov-coalgebra}
Let $(C,\top)$ be a locality vector space.
A {\bf locality Novikov coproduct} is a linear map
$\delta:C\to C\tloc C$ satisfying the polar decomposition property
\[
 \delta(U^\top)\subseteq U^\top\tloc U^\top
 \qquad(U\subseteq C),
\]
and
\begin{align*}
 (\delta\otimes\id)\delta-(\id\otimes\delta)\delta
 &=\flip_{12}\bigl((\delta\otimes\id)\delta
                        -(\id\otimes\delta)\delta\bigr),\\
 (\delta\otimes\id)\delta
 &=\flip_{23}(\delta\otimes\id)\delta.
\end{align*}
These are \textup{(NC1)} and \textup{(NC2)} as identities
$C\to C^{\otimes_\top3}$.

A {\bf locality multi-Novikov coalgebra} is a tuple
$(C,\top,(\delta_\omega)_{\omega\in\Omega})$, where each
$\delta_\omega:C\to C\tloc C$ is linear and satisfies
\begin{equation}\mlabel{eq:locality-multi-co-polar}
 \delta_\omega(U^\top)\subseteq U^\top\tloc U^\top
 \qquad(U\subseteq C,\ \omega\in\Omega).
\end{equation}
For every $\omega,\tau\in\Omega$, put
\[
 \mathfrak a_{\omega,\tau}
 =(\delta_\omega\otimes\id)\delta_\tau
   -(\id\otimes\delta_\tau)\delta_\omega,
\]
The required identities $C\to C^{\otimes_\top3}$ are
\begin{align*}
 \mathfrak a_{\omega,\tau}
   &=\flip_{12}\mathfrak a_{\omega,\tau},
       &&\tag{MNC1}\\
 \mathfrak a_{\omega,\tau}
   &=\mathfrak a_{\tau,\omega},
       &&\tag{MNC2}\\
 (\delta_\omega\otimes\id)\delta_\tau
   &=\flip_{23}(\delta_\tau\otimes\id)\delta_\omega.
       &&\tag{MNC3}
\end{align*}
Omitting \textup{(MNC2)} gives a {\bf noncommuting locality
multi-Novikov coalgebra}.
\end{definition}

\begin{remark}
The polar decomposition property gives, for every $U\subseteq C$,
\[
\begin{aligned}
 (\delta_\omega\otimes\id)\delta_\tau(U^\top)
       &\subseteq (U^\top)^{\otimes_\top3},\\
 (\id\otimes\delta_\tau)\delta_\omega(U^\top)
       &\subseteq (U^\top)^{\otimes_\top3},
\end{aligned}
\]
\end{remark}

\begin{proposition}
\mlabel{pp:locality-co-Gelfand}
Let $(C,\top,\Delta)$ be a coassociative locality coalgebra, not
necessarily counital, and let $(D_\omega)_{\omega\in\Omega}$ be a family of locality coderivations. Set
\[
 \delta_\omega=(\id\otimes D_\omega)\Delta,
\]
Then each $\delta_\omega$ has the polar decomposition property:
\[
 \delta_\omega(U^\top)\subseteq U^\top\tloc U^\top
 \qquad(U\subseteq C).
\]
For all $\omega,\tau\in\Omega$, the mixed Novikov defects are
\begin{align}
 \mathfrak a_{\omega,\tau}
       -\flip_{12}\mathfrak a_{\omega,\tau}
 &=-(\id\otimes\id\otimes D_\tau D_\omega)
          (\id-\flip_{12})\Delta^{(2)},
 \mlabel{eq:locality-co-cocommutativity-defect-one}\\
 \mathfrak a_{\omega,\tau}-\mathfrak a_{\tau,\omega}
 &=(\id\otimes\id\otimes[D_\omega,D_\tau])\Delta^{(2)},
 \mlabel{eq:locality-co-commutator-defect}\\
 (\delta_\omega\otimes\id)\delta_\tau
       -\flip_{23}(\delta_\tau\otimes\id)\delta_\omega
 &=(\id\otimes D_\omega\otimes D_\tau)
          (\id-\flip_{23})\Delta^{(2)}.
 \mlabel{eq:locality-co-cocommutativity-defect-three}
\end{align}
For every $U\subseteq C$, both sides of these identities map
$U^\top$ into $(U^\top)^{\otimes_\top3}$.
\end{proposition}
\begin{proof}
Fix $U\subseteq C$ and $c\in U^\top$. Choose
\[
 \Delta c=\sum_j a_j\otimes b_j,\qquad
 a_j,b_j\in U^\top,\quad a_j\top b_j.
\]
Polar stability gives
\[
 D_\omega b_j\in U^\top\cap\{a_j\}^\top,\qquad
 \delta_\omega c
 =\sum_j a_j\otimes D_\omega b_j\in U^\top\tloc U^\top.
\]
For each $j$,
\[
 a_j\in(U\cup\{b_j\})^\top
 \quad\Longrightarrow\quad
 \Delta a_j\in
 (U\cup\{b_j\})^\top\tloc(U\cup\{b_j\})^\top.
\]
Thus
\[
 \Delta^{(2)}c=\sum c_1\otimes c_2\otimes c_3
       \in(U^\top)^{\otimes_\top3},
 \qquad c_i\top c_j\quad(i\ne j)
\]
in every summand. Polar stability of finite compositions gives
\[
 \begin{gathered}
  c_1\otimes D_\omega c_2\otimes D_\tau c_3
       \in(U^\top)^{\otimes_\top3},\\
  c_1\otimes c_2\otimes D_\tau D_\omega c_3
       \in(U^\top)^{\otimes_\top3},\qquad
  c_1\otimes c_2\otimes D_\omega D_\tau c_3
       \in(U^\top)^{\otimes_\top3}.
 \end{gathered}
\]
The coderivation identity therefore yields, on these domains,
\begin{align*}
 (\delta_\omega\otimes\id)\delta_\tau(c)
 &=\sum c_1\otimes D_\omega c_2\otimes D_\tau c_3,\\
 (\id\otimes\delta_\tau)\delta_\omega(c)
 &=\sum c_1\otimes D_\omega c_2\otimes D_\tau c_3
   +\sum c_1\otimes c_2\otimes D_\tau D_\omega c_3,\\
 \mathfrak a_{\omega,\tau}(c)
 &=-\sum c_1\otimes c_2\otimes D_\tau D_\omega c_3.
\end{align*}
Moreover,
\[
 \flip_{12}\bigl((U^\top)^{\otimes_\top3}\bigr)
 =(U^\top)^{\otimes_\top3}
 =\flip_{23}\bigl((U^\top)^{\otimes_\top3}\bigr).
\]
Proposition~\mref{pp:ordinary-co-Gelfand} now gives the three
displayed defect formulas, with every term in the stated
restricted tensor space.
\end{proof}

\begin{corollary} If $(C, \top, \Delta)$ is cocommutative, then $(C, \top ,\Delta, \delta _\omega)$ is a noncommuting multi-Novikov
coalgebra. In this case, if the coderivations are pairwise commuting, then $(C, \top , \Delta, \delta _\omega)$ is a  multi-Novikov coalgebra.
\end{corollary}

\begin{corollary} If $(C, \top, \Delta)$
has a counit, then, without assuming the cocommutativity, we have
\[
 \mathfrak a_{\omega,\tau}=\mathfrak a_{\tau,\omega}
 \quad\Longleftrightarrow\quad [D_\omega,D_\tau]=0.
\]
\end{corollary}

\section{Induced locality Novikov structures for bialgebras}
\mlabel{sec:locality-novikov}
The principal structure of this section is the one induced by a
locality biderivation of a {\bf commutative} locality bialgebra.

\subsection{Locality bialgebras with Novikov pairs}
\mlabel{subsec:Novikov-pairs}
\begin{definition}\mlabel{def:locality-bialgebra}
A {\bf locality bialgebra} \cite[Definition~5.1]{CGPZRen}
is a tuple $(B,\top,m,u,\Delta,\eps)$, where
$(B,\top,m,u)$ is a unital associative locality algebra,
$(B,\top,\Delta,\eps)$ is a counital coassociative locality coalgebra, and
$\Delta$ and $\eps$ are locality algebra homomorphisms, that is,
\begin{equation}\mlabel{eq:bialgebra-multiplicativity}
 \begin{gathered}
  \Delta(1)=1\otimes1,\qquad \eps(1)=1,\\
  \Delta(xy)=\Delta(x)\sqprod\Delta(y),\qquad
  \eps(xy)=\eps(x)\eps(y)\qquad(x\top y).
 \end{gathered}
\end{equation}
where $1=u(1_\kk)$.

It is {\bf commutative} if
\[
 m(x,y)=m(y,x)\qquad(x\top y).
\]
\end{definition}

\begin{definition}\mlabel{def:hopf-bider}
Let $(B,\top,m,u,\Delta,\eps)$ be a locality bialgebra.
A linear map $D:B\to B$ is a {\bf locality biderivation} if
\[
\begin{gathered}
 D(U^\top)\subseteq U^\top\qquad(U\subseteq B),\\
 D(xy)=D(x)y+xD(y)\qquad(x\top y),\\
 \Delta D=(D\otimes\id+\id\otimes D)\Delta,
 \qquad \eps D=0.
\end{gathered}
\]
A tuple
\[
 (B,\top,m,u,\Delta,\eps,(D_\omega)_{\omega\in\Omega})
\]
is a {\bf commuting locality multidifferential bialgebra}
if each $D_\omega$ is a locality biderivation and
\[
 [D_\omega,D_\tau]=0\qquad(\omega,\tau\in\Omega).
\]
\end{definition}

\begin{definition}\mlabel{def:locality-Novikov-pair}
A {\bf commutative locality bialgebra with a Novikov pair} is a tuple
$(B,\top,m,\Delta,\eps,\tri,\rho)$ satisfying the following conditions.
\begin{enumerate}[label=\textup{(B\arabic*)}]
\item \mlabel{i:b1} The sextuple $(B,\top,m,u,\Delta,\eps)$ is a commutative locality bialgebra.
\item \mlabel{i:b2} The product $\tri:B\times_\top B\to B$ is locality bilinear
and satisfies the {\bf polar stability}
\[
 a,b\in U^\top,\ a\top b
 \quad\Longrightarrow\quad a\tri b\in U^\top
 \qquad(U\subseteq B).
\]
\item \mlabel{i:b3} For pairwise independent $a,b,c$, we have
\begin{align}
 (ab)\tri c&=a(b\tri c),
 \mlabel{eq:based-locality-product-one}\\
 a\tri(bc)&=(a\tri b)c+b(a\tri c).
 \mlabel{eq:based-locality-product-two}
\end{align}
\item \mlabel{i:b4} The map $\rho:B\to B\tloc B$ is linear and satisfies
the {\bf polar decomposition}
\begin{equation}\mlabel{eq:based-locality-polar}
 \rho(U^\top)\subseteq U^\top\tloc U^\top
 \qquad(U\subseteq B).
\end{equation}
\item \mlabel{i:b5} The following identities hold in $B^{\otimes3}:$
\begin{align}
 (\Delta\otimes\id)\rho&=(\id\otimes\rho)\Delta,
 \mlabel{eq:based-locality-co-one}\\
 (\id\otimes\Delta)\rho
 &=(\rho\otimes\id)\Delta+(\id\otimes\rho)\Delta.
\mlabel{eq:based-locality-co-two}
\end{align}
\item \mlabel{i:b6} There is the {\bf unit-counit compatibility}.
\begin{equation}\mlabel{eq:based-locality-sync}
 1\tri a=(\eps\otimes\id)\rho(a)\qquad(a\in B).
\end{equation}
\end{enumerate}
The pair is denoted by $(\tri,\rho)$.
\end{definition}

\begin{proposition} The product $\tri$ in Definition~\mref{def:locality-Novikov-pair} is a locality Novikov product.
\end{proposition}

\begin{proof}
Let $D:=1\tri(-)$. Since $\{1\}^\top=B$, the locality bilinearity and
\mref{i:b2} imply that $D$ is in $\mathrm{End}_{\kk}(B)$ and $ D(U^\top)\subseteq U^\top$ for $U\subseteq B$.
If $a\top b$, then taking the triples $(a,1,b)$ and $(1,a,b)$ in
\meqref{eq:based-locality-product-one} and \meqref{eq:based-locality-product-two}
give
\[
 a\tri b=a(1\tri b)=aD(b),\qquad
 D(ab)=D(a)b+aD(b).
\]
For pairwise independent $a,b,c$, the polar stability makes all products below defined, and the commutativity and locality associativity establish the equality
\[
 (a\tri b)\tri c=aD(b)D(c)
 =aD(c)D(b)=(a\tri c)\tri b.
\]
For the associator,
\[
 \begin{aligned}
 A_\tri(a,b,c)
 &:=(a\tri b)\tri c-a\tri(b\tri c)\\
 &=aD(b)D(c)-aD\bigl(bD(c)\bigr)\\
 &=-abD^2(c)=-baD^2(c)=A_\tri(b,a,c).
 \end{aligned}
\]
Together with the locality bilinearity and polar stability in
\mref{i:b2}, these are the locality Novikov axioms.
\end{proof}

\begin{definition}\mlabel{def:Novikov-pair-morphism}
Let $f:B\to B'$ be a linear map between locality bialgebras. Then $f$ is

\begin{enumerate}
    \item A {\bf locality bialgebra homomorphism} if
     \[
    a\top b\Longrightarrow f(a)\top' f(b),\qquad f(1)=1',
    \]
     \[
     f(ab)=f(a)f(b)\quad(a\top b),\qquad
     \Delta'f=(f\otimes f)\Delta,\qquad\eps'f=\eps.
    \]

    \item A {\bf differential locality bialgebra homomorphism} if additionally $fD=D'f$.

    \item A {\bf morphism of Novikov pairs} if additionally
    \[
    f(a\tri b)=f(a)\tri' f(b)\quad(a\top b),\qquad
    (f\otimes f)\rho=\rho'f.
    \]

\end{enumerate}
\end{definition}

\begin{theorem}
\mlabel{thm:bialgebra-Novikov-recovery}
Let $(B,\top,m,u,\Delta,\eps)$ be a commutative locality bialgebra. A locality biderivation and a locality Novikov pair can be mutually determined as follows:
\begin{enumerate}[label=\textup{(\roman*)}]
\item For a locality biderivation $D$, set
\begin{equation}\mlabel{eq:based-locality-induction}
 a\tri_D b:=aD(b)\quad(a\top b),\qquad
 \rho_D:=(\id\otimes D)\Delta.
\end{equation}
\item For a locality Novikov pair $(\tri,\rho)$, set
\begin{equation}\mlabel{eq:based-locality-recovery}
 D_\tri(a):=1\tri a=(\eps\otimes\id)\rho(a).
\end{equation}
\end{enumerate}
The constructions also identify morphisms of Novikov pairs with
differential locality bialgebra homomorphisms.
\end{theorem}
\begin{proof} Let $D$ be a biderivation. Polar stability of $D$ gives
\[
 a\top b\Longrightarrow a\top D(b),\qquad
 a,b\in U^\top\Longrightarrow a,D(b)\in U^\top.
\]
Hence $a\tri_D b=aD(b)$ is polar stable.

For pairwise independent $a,b,c$, the locality associativity and locality commutativity give
\[
 (ab)\tri_D c=abD(c)=a(b\tri_D c),
\]
\[
 a\tri_D(bc)=aD(b)c+abD(c)
             =(a\tri_D b)c+b(a\tri_D c).
\]

Fix $U\subseteq B$ and $h\in U^\top$, the polar decomposition property of
$\Delta$ gives
\[
 \Delta h=\sum_j a_j\otimes b_j, \qquad
 a_j,b_j\in U^\top,\quad a_j\top b_j.
\]
Then
\[
 D(b_j)\in U^\top\cap\{a_j\}^\top,
 \qquad
 \rho_D(h)=\sum_j a_j\otimes D(b_j)\in U^\top\tloc U^\top.
\]
Thus $\rho_D$ has the polar decomposition property.

By Lemma~\mref{lem:coproduct-square-locality}, write $\Delta^{(2)}h=\sum h_1\otimes h_2\otimes h_3$ with pairwise independent factors. Then
\begin{align*}
 (\Delta\otimes\id)\rho_D(h)
 &=\sum h_1\otimes h_2\otimes D(h_3)
  =(\id\otimes\rho_D)\Delta(h),\\
 (\id\otimes\Delta)\rho_D(h)
 &=\sum h_1\otimes D(h_2)\otimes h_3
   +\sum h_1\otimes h_2\otimes D(h_3)\\
 &=\bigl((\rho_D\otimes\id)\Delta
          +(\id\otimes\rho_D)\Delta\bigr)(h).
\end{align*}

Finally,
\[
 (\eps\otimes\id)\rho_D
 =D(\eps\otimes\id)\Delta=D=1\tri_D(-),
\]
which proves synchronization.

Now let $(\tri, \rho)$ be a Novikov pair. Set $D_A(a)=1\tri a$ and $D_C(a)=(\eps\otimes\id)\rho(a)$.

Letting $b=1$ in  (\mref {eq:based-locality-product-one}), and letting $a=1$ in  (\mref {eq:based-locality-product-two}), we have
\[
 a\tri c=a(1\tri c)=aD_A(c),
\]
\[
 D_A(bc)=(1\tri b)c+b(1\tri c)
          =D_A(b)c+bD_A(c)\qquad(b\top c).
\]
For $a\in U^\top$, we have $1,a\in U^\top$ and $1\top a$.
Consequently, $D_A(a)=1\tri a\in U^\top$, and so $D_A$ is a locality derivation. The Leibniz rule gives
$ D_A(1)=2D_A(1)$, hence $D_A(1)=0$.

Next, apply $\id\otimes\eps\otimes\id$ to
\meqref{eq:based-locality-co-one}. Writing
$\rho(h):=\sum r_j\otimes s_j$ and
$\Delta h:=\sum h_1\otimes h_2$, then
\begin{align*}
 (\id\otimes\eps\otimes\id)(\Delta\otimes\id)\rho(h)
 &=\sum r_j\otimes s_j=\rho(h),\\
 (\id\otimes\eps\otimes\id)(\id\otimes\rho)\Delta(h)
 &=\sum h_1\otimes D_C(h_2).
\end{align*}
Thus
\begin{equation}\mlabel{eq:recovered-rho}
 \rho=(\id\otimes D_C)\Delta.
\end{equation}

Apply $\eps\otimes\id\otimes\id$ to \meqref{eq:based-locality-co-two}. The three terms become
\[
 \Delta D_C,\qquad (D_C\otimes\id)\Delta,\qquad\rho.
\]
Thanks to \meqref{eq:recovered-rho}, we have
\[
 \Delta D_C=(D_C\otimes\id+\id\otimes D_C)\Delta.
\]
For $q=\eps D_C$, the counitality implies
$
 q=(\eps\otimes\eps)\Delta D_C=q+q$ and so $q=0$.
For $h\in U^\top$, the polar decomposition \mref{i:b4} yields
\[
 \rho(h)=\sum_j r_j\otimes s_j,\qquad r_j,s_j\in U^\top,
 \qquad
 D_C(h)=\sum_j\eps(r_j)s_j\in U^\top.
\]
Thus $D_C$ is a locality coderivation. The unit-counit compatibility \mref{i:b6} gives
$D_A=D_C=:D$.

Thus we have completed the proof that biderivations and Novikov pairs are mutually determined.

For a locality bialgebra map $f:B\to B'$, preservation of $\tri$
implies
\[
 fD(a)=f(1\tri a)=1'\tri'f(a)=D'f(a).
\]
Conversely, if $fD=D'f$, then, for $a\top b$,
\[
 f(a\tri b)=f(aD(b))=f(a)D'f(b)=f(a)\tri'f(b),
\]
\[
 (f\otimes f)\rho
 =(\id\otimes D')(f\otimes f)\Delta
 =(\id\otimes D')\Delta'f=\rho'f.
\]
The target products are defined because $f$ preserves independence.
Preservation of $\rho$ alone also gives
\[
 D'f=(\eps'\otimes\id)\rho'f
     =(\eps'\otimes\id)(f\otimes f)\rho=fD.
\]
Both constructions retain the same base maps. They are inverse on
objects and morphisms.
\end{proof}

\subsection{Locality Hopf algebras with Novikov pairs}
\mlabel{subsec:based-novikov-hopf}
For the Hopf construction, we add a tensor-factorization condition for the locality multiplication so that the convolution and antipode contractions are well-defined linear maps, because our coproduct has the polar decomposition property.

\begin{definition}
\mlabel{def:locality-Hopf}
\begin{enumerate}
    \item
If a locality bilinear multiplication
$$m: H\times _\top H\to H$$
factorizes through
\begin{equation}\mlabel{eq:restricted-multiplication}
 \mu:H\tloc H\longrightarrow H,\qquad
 \mu(x\otimes y)=m(x,y)\quad(x\top y),
\end{equation}
then $m$ is said to have the {\bf tensor-factorization property.}
\item
A {\bf locality Hopf algebra} is a locality bialgebra $(H,\top,m,u,\Delta,\eps)$ whose multiplication $m$ has the tensor-factorization property and which is endowed with an antipode $S$ in the sense of
\cite[Definition~5.3]{CGPZRen}, that is, the locality linear map $S:H\to H$ satisfies
\[
 x\top y\ \Longrightarrow\ S(x)\top y
 \qquad(x,y\in H),
\]
\begin{equation}\mlabel{eq:antipode-identities}
 \mu(S\otimes\id)\Delta=u\eps=\mu(\id\otimes S)\Delta.
\end{equation}
\item
For a commuting family of locality biderivations
$(D_\omega)_{\omega\in\Omega}$, the tuple
\[
 (H,\top,m,u,\Delta,\eps,S,(D_\omega)_{\omega\in\Omega})
\]
is called a {\bf commuting locality multidifferential Hopf algebra}.
\end{enumerate}
\end{definition}

\begin{remark}
For polar-stable linear maps $f,g:H\to H$, one has
\[
 (f\otimes g)\Delta(H)\subseteq H\tloc H,
\]
The multiplication factorization allows the product
\[
 f*g=\mu(f\otimes g)\Delta,
\]
to be well-defined independent of tensor expansions.
\end{remark}

\begin{definition}
\mlabel{def:based-locality-Novikov-Hopf}
A {\bf locality Hopf algebra with a Novikov pair} is a tuple
\[
 (H,\top,m, u,\Delta,\eps,S,\tri,\rho)
\]
whose base $(H,\top,m,u,\Delta,\eps,S)$ is a locality Hopf algebra
as in Definition~\mref{def:locality-Hopf}, and whose pair
$(\tri,\rho)$ satisfies Definition~\mref{def:locality-Novikov-pair}.
\end{definition}

\begin{definition}
\mlabel{def:based-locality-Novikov-Hopf-morphism}
A {\bf locality Hopf map} is a locality bialgebra map between
locality Hopf algebras that intertwines their antipodes.
A morphism of locality Hopf algebras with Novikov pairs is a morphism of
the underlying Novikov pairs, as in
Definition~\mref{def:Novikov-pair-morphism}, that also satisfies
\[
 fS=S'f.
\]
The antipode condition is automatic for these locality bialgebra
maps, as proved below.
\end{definition}

\begin{proposition}\mlabel{prop:antipode}
Let $H$ be a locality Hopf algebra, and $D$ be a locality biderivation. Then $DS=SD$.
\end{proposition}
\begin{proof}
Applying $D$ to $\mu(S\otimes\id)\Delta=u\eps$ and using $D(1)=0$ gives
\[
 \sum D(S(x_{(1)}))x_{(2)}+\sum S(x_{(1)})D(x_{(2)})=0.
\]
On the other hand, applying the antipode identity to $D(x)$ and using the coderivation identity, we obtain
\[
 \sum S(D(x_{(1)}))x_{(2)}+\sum S(x_{(1)})D(x_{(2)})=0.
\]
Let
$$F:=DS-SD.$$
A direct check gives
$F*\id=0$.
Consequently,
$$(F*\id)*S=F*(\id*S).$$
Hence
$0=F*(u\eps)=F$, as required.
\end{proof}

Since no additional compatibility condition was imposed between the antipode and biderivation or between the antipode and Novikov pair, by Theorem~\mref{thm:bialgebra-Novikov-recovery}, we have

\begin{corollary}
\mlabel{thm:based-locality-recovery}
Let $H=(H,\top,m,\Delta,\eps,S)$ be a commutative locality Hopf
algebra. Let
$$ \mathcal{D}:=\{(H,D)\,|\, D \text{ is a biderivation on } H\},$$
$$\mathcal{P}:= \{(H,\tri,\rho)\,|\, (\tri,\rho) \text{ is a Novikov pair on } H\}.$$
Then  the assignments
\[
\begin{gathered}
 D\longmapsto(\tri_D,\rho_D),\qquad
 a\tri_D b:=aD(b)\ (a\top b),\quad
 \rho_D:=(\id\otimes D)\Delta,\\
 (\tri,\rho)\longmapsto D,\qquad
 D(a):=1\tri a=(\eps\otimes\id)\rho(a),
\end{gathered}
\]
define a one-to-one correspondence between $\mathcal{D}$ and $\mathcal{P}$.
\end{corollary}

\section{Chen fractions and the MZV shuffle Hopf algebra}\mlabel{sec:mzv}\mlabel{sec:chen-hopf}
We now apply the locality constructions to Chen fractions.
Their intrinsic variable supports control independence, while partial
derivatives raise the denominator exponents without enlarging support.
This produces explicit locality Novikov and Hopf structures, and finally gives algebraic structures for MZVs.

\subsection{Chen fractions}\mlabel{subsec:positive-chen-derived}
For pairwise distinct positive integer labels $\bfi=(i_1,\ldots,i_k)$ and
$\mathbf s=(s_1,\ldots,s_k)\in\Z_{\ge1}^k$, we define the Chen fraction
\begin{equation}\mlabel{eq:chen}
 F_{\mathbf s;\bfi}=\frac 1{(x_{i_1}+\cdots+x_{i_k})^{s_1}(x_{i_2}+\cdots+x_{i_k})^{s_2}\cdots x_{i_k}^{s_k}}.
\end{equation}
Let $F^\Ch$ be the set consisting of these fractions and $1$, and let $\calf^\Ch$ be the $\Q$-linear subspace of functions consisting of the $\Q$-linear combinations of these fractions.

\begin{proposition} \cite {GPZGalois} $F^\Ch$ is a $\Q$-linear basis of $\calf ^\Ch$.
\end{proposition}

We define
$$\rmsupp (F_{\mathbf s;\bfi} )=\{i_1, \cdots, i_k\}$$
and
$$\rmsupp (\sum a_{\mathbf s;\bfi}F_{\mathbf s;\bfi} )=\bigcup _{a_{\mathbf s;\bfi}\ne 0}\rmsupp (F_{\mathbf s;\bfi}) $$
with $\rmsupp (0)= \rmsupp (1)=\varnothing$, which is well-defined because $F^\Ch$ is a basis.

Then in $\calf ^\Ch$, define
\begin{equation}\mlabel{eq:chenloc}
 f\top g\quad\Longleftrightarrow\quad\rmsupp (f)\cap\rmsupp (g)=\varnothing.
\end{equation}

\begin{proposition}
 \mcite{GHXZ1,GHXZ,GX} $(\calf  ^\Ch, \top, \cdot)$ is a commutative associative locality algebra.
\end{proposition}

Multiplication in the ambient rational-function field is bilinear.
Since products of independent Chen fractions remain in $\calf^\Ch$,
it induces a linear map
\[
 \mu:\calf^\Ch\tloc\calf^\Ch\longrightarrow\calf^\Ch,
 \qquad \mu(f\otimes g)=fg\quad(f\top g).
\]
Thus the additional factorization~\meqref{eq:restricted-multiplication}
holds for the Chen fractions.

Put $\partial_j=-\partial/\partial x_j$ for $j\ge1$, and henceforth set $\partial_0=0$.
If $j$ is absent from the label set, then $\partial_jF_{\mathbf s;\bfi}=0$. If $j=i_q$, then
\begin{equation}\mlabel{eq:chen-der}
 \partial_{i_q}F_{\mathbf s;\bfi}
 =\sum_{r=1}^q s_rF_{\mathbf s+\mathbf e_r;\bfi}.
\end{equation}
If we put $d_{a,b}=\partial_a-\partial_b$, then
\begin{equation}\mlabel{eq:d-raise-local}
 d_{i_q,i_{q-1}}F_{\mathbf s;\bfi}=s_qF_{\mathbf s+\mathbf e_q;\bfi}.
\end{equation}

Consequently,
\begin{equation}\mlabel{eq:chen-support-derivative}
 \partial_j(\calf ^\Ch)\subseteq\calf ^\Ch,\qquad
 \rmsupp (\partial_j f)\subseteq\rmsupp (f).
\end{equation}
For any scalar sequence $c=(c_j)_{j\ge1}\in\prod_{j\ge1}\Q$, define the pointwise finite sum
\begin{equation}\mlabel{eq:weighted-chen-derivation}
 D_cf=\sum_{j\ge1}c_j\partial_jf.
\end{equation}
This is a locality derivation, and $D_cD_d=D_dD_c$.

\begin{corollary}\mlabel{thm:chen-MN}
For $f\top g$, define $f\tri_cg=fD_c(g)$. Then
$(\calf ^\Ch,\top,\{\tri_c\}_{c\in\prod\Q})$ is a locality multi-Novikov algebra.
\end{corollary}

Define the coproduct \mcite{GHXZ} on all-one fractions by deconcatenation:
\begin{equation}\mlabel{eq:chen-coprod-one}
 \Delta_\Ch F_{\mathbf1;(i_1,\ldots,i_k)}
 =\sum_{q=0}^k F_{\mathbf1;(i_1,\ldots,i_q)}\otimes
 F_{\mathbf1;(i_{q+1},\ldots,i_k)}.
\end{equation}
and extend it recursively by
\begin{equation}\mlabel{eq:chen-coprod-recursion}
 \Delta_\Ch F_{\mathbf s+\mathbf e_q;\bfi}
 =\frac1{s_q}\bigg(d_{i_q,i_{q-1}} \otimes\id+\id\otimes d_{i_q,i_{q-1}}\bigg)\Delta_\Ch F_{\mathbf s;\bfi}.
\end{equation}
For the fixed label sequence $\mathbf i$, put
\[
T_q:=d_{i_q,i_{q-1}}\otimes\id
+\id\otimes d_{i_q,i_{q-1}},\qquad 1\le q\le k.
\]
The operators $T_q$ commute, since the partial derivations commute.
Thus the recursion has the order-independent form
\[
\Delta_\Ch F_{\mathbf s;\mathbf i}
=\left(\prod_{q=1}^k
\frac{T_q^{s_q-1}}{(s_q-1)!}\right)
\Delta_\Ch F_{\mathbf1;\mathbf i}.
\]
The fact that these formulas define a locality Hopf coproduct is
part of the Chen-fraction Hopf construction cited below.
The counit takes the value $1$ on the unit and vanishes on nonconstant Chen fractions.

\begin{proposition}
[\mcite{GHXZ}]\mlabel{thm:known-chen-hopf}
The coproduct and counit above make $\calf ^\Ch$ a connected, weight-graded locality commutative Hopf algebra.
Every $\partial_j$ is a locality biderivation. In particular,
\begin{equation}\mlabel{eq:chen-coderivation}
 \Delta_\Ch\partial_j=(\partial_j\otimes\id+\id\otimes\partial_j)\Delta_\Ch.
\end{equation}
\end{proposition}

\begin{remark}
The coproduct of a depth-two all-one fraction contains the ordered middle term
$F_{\mathbf1;(i_1)}\otimes F_{\mathbf1;(i_2)}$, without its transposed term.
Thus $\Delta_\Ch$ is generally not cocommutative.
\end{remark}

\begin{corollary}
\mlabel{cor:chen-co-derived}
For
$D=D_c$ of \meqref{eq:weighted-chen-derivation}, put
\[
 f\tri_c g=fD_c(g)\quad(f\top g),\qquad
 \delta_c =(\id\otimes D_c)\Delta_\Ch.
\]
Then $(\calf^\Ch,\top,m,\Delta_\Ch,\eps,S,
\tri_c,\delta_c)$ is a locality Hopf algebra with a Novikov pair
in the sense of Definition~\mref{def:based-locality-Novikov-Hopf}.
\end{corollary}

\subsection{Differential MZV shuffle Hopf algebra}\mlabel{subsec:mzv-novikov-descent}
We end the paper by giving the application of the locality framework developed in this paper to a non-local setting, yielding a biderivation on the shuffle Hopf algebra of multiple zeta values. Further details will be provided separately.

The locality of labelled fractions is essential for their multiplication.
%After forgetting the labels, the shuffle product becomes everywhere defined; only certain combinations of raising operators are compatible with this forgetting map.
Let
\[
 \mathcal H_{\Z_{\ge1}}=\Q1\oplus
 \bigoplus_{k\ge1,\,\mathbf s\in\Z_{\ge1}^k}\Q[\mathbf s]
\]
equipped with the MZV shuffle product $\shap$.
Forgetting labels gives a surjective locality algebra homomorphism
\[
 \pi:\calf ^\Ch\to\mathcal H_{\Z_{\ge1}},\qquad
 F_{\mathbf s;\bfi}\longmapsto[\mathbf s],\quad\pi(1)=1.
\]
The Chen fraction coproduct descends to the coproduct $\Delta_{\rm sh}$ of \mcite{GHXZ}, satisfying
$\Delta_{\rm sh}\pi=(\pi\otimes\pi)\Delta_\Ch$.

\begin{proposition}\mlabel{prop:symmetric-chen-descent}
Denote $\nabla_{\rm tot}:=\sum_{j\ge1}\partial_j$. Define a linear operator $\mathscr{R}$ on $\mathcal H_{\Z_{\ge1}}$ by
\begin{equation}\mlabel{eq:mzv-total-raising}
 \mathscr R(1):=0,\qquad
 \mathscr R[\mathbf s]:=\sum_{r=1}^k(k-r+1)s_r[\mathbf s+\mathbf e_r], \quad [\mathbf s]=[s_1, \cdots, s_k].
\end{equation}
Then $\mathscr R$ is the descent of $\nabla_{\rm tot}$:
$$\pi\nabla_{\rm tot}=\mathscr R\pi.$$
Furthermore, $\mathscr R$ is a derivation and coderivation:
\begin{align}
 \mathscr R(u\shap v)&=\mathscr R(u)\shap v+u\shap\mathscr R(v),\mlabel{eq:R-shuffle-derivation}\\
 \Delta_{\rm sh}\mathscr R
 &=(\mathscr R\otimes\id+\id\otimes\mathscr R)\Delta_{\rm sh}.\mlabel{eq:R-shuffle-coderivation}
\end{align}
\end{proposition}

\begin{proof}
By \meqref{eq:chen-der}, we have
\[
 \nabla_{\rm tot}F_{\mathbf s;\bfi}
 =\sum_{q=1}^k\sum_{r=1}^q s_rF_{\mathbf s+\mathbf e_r;\bfi}
 =\sum_{r=1}^k(k-r+1)s_rF_{\mathbf s+\mathbf e_r;\bfi}.
\]
Note that the coefficients do not depend on the names of the labels. Thus applying $\pi$ gives the descent formula.

For basis elements $u, v\in \mathcal H_{\Z_{\ge 1}}$, choose their lifts $f,g$ with disjoint label sets. The Leibniz rule and locality multiplicativity of $\pi$ give
\[
 \mathscr R(u\shap v)=\pi\nabla_{\rm tot}(fg)
 =\pi\bigl(\nabla_{\rm tot}(f)g+f\nabla_{\rm tot}(g)\bigr),
\]
which proves \meqref{eq:R-shuffle-derivation}.
Also, each $\partial_j$ is a coderivation, and every element of $\calf ^\Ch$ involves only finitely many variables.
Thus the pointwise finite sum $\nabla_{\rm tot}$ is also a coderivation.
Applying $\pi\otimes\pi$ and using surjectivity of $\pi$ proves \meqref{eq:R-shuffle-coderivation}.
\end{proof}

\begin{example}
For $s,a,b,c\in\Z_{\ge1}$, formula~\meqref{eq:mzv-total-raising} gives
\begin{align*}
 \mathscr R[s]&=s[s+1],\quad
 \mathscr R[a,b]=2a[a+1,b]+b[a,b+1],\\
 \mathscr R[a,b,c]&=3a[a+1,b,c]+2b[a,b+1,c]+c[a,b,c+1].
\end{align*}
\end{example}

\begin{corollary}\mlabel{cor:mzv-shuffle-novikov}
Define
\begin{equation}\mlabel{eq:mzv-shuffle-novikov}
 u\tri _{\rm sh}v:=u\shap\mathscr R(v).
\end{equation}
\begin{enumerate}
    \item
Then $\mathcal H_{\Z_{\ge1}}$ is a Novikov algebra, and
\begin{align}
 \Delta_{\rm sh}(u\tri _{\rm sh}v)
 ={}&\sum(u_{(1)}\tri _{\rm sh}v_{(1)})\otimes(u_{(2)}\shap v_{(2)})
+\sum(u_{(1)}\shap v_{(1)})\otimes(u_{(2)}\tri _{\rm sh}v_{(2)}).
 \mlabel{eq:mzv-novikov-hopf-compatibility}
\end{align}
\item
If $\mathcal H_{w,d}$ denotes the homogeneous subspace of weight $w$ and depth $d$, then
\[
 \mathcal H_{w,d}\tri _{\rm sh}\mathcal H_{w',d'}
 \subseteq\mathcal H_{w+w'+1,d+d'}.
\]
\item Let $\mathbf{c}=(c_j)_{j\ge 1}\in \prod_{j\ge 1}\Q$ with $c_j=1$ for $j\ge 1$, so that $ \nabla_{\rm tot}=D_{\mathbf{c}}$. Then
$$ \tri_{\rm sh} (\pi\otimes \pi)  = \pi\, \tri_{\mathbf{c}}.$$
\end{enumerate}
\end{corollary}

Since a Novikov algebra is in particular a pre-Lie algebra, we obtain
\begin{corollary}\mlabel{cor:mzv-novikov-lie}
The product
\begin{equation}\mlabel{eq:mzv-novikov-lie}
 [u,v]_{\mathscr R}:=u\shap\mathscr R(v)-v\shap\mathscr R(u)
\end{equation}
defines a Lie bracket on $\mathcal H_{\Z_{\ge1}}$, with $[1,[\mathbf s]]_{\mathscr R}=\mathscr R[\mathbf s]$.
\end{corollary}

\noindent
{\bf Acknowledgments.} This research is supported by NSFC (12471062).

\noindent
{\bf Declaration of interests.} The authors have no conflict of interest to declare.

\noindent
{\bf Data availability.} No data were created or analyzed.

\end{document}